\documentclass[11pt]{article}
\usepackage{amssymb,amsmath}
\usepackage[all]{xy}

\usepackage{graphicx}
\begin{document}

\title{Symbolic Dynamics Generated by a Combination of Graphs}
\author{Vasileios Basios$^1$, Gian-Luigi Forti$^{1,2}$ and Gregoire Nicolis$^1$}
\maketitle
\begin{center}
$^1$ Interdisciplinary Center for Nonlinear Phenomena and
Complex Systems,
C.P. 231, Universit\'e  Libre de Bruxelles, B-1050 Brussels, Belgium.\\

Contact: \tt{vbasios@ulb.ac.be}, \tt{gnicolis@ulb.ac.be}
\end{center}
\begin{center}
$^2$ Dipartimento di Matematica, Universit\`a degli Studi di Milano,\\
 via C. Saldini 50, I-20133,
Milano, Italy.
Contact: \tt{forti@mat.unimi.it}
\end{center}


\begin{abstract}

In this paper we investigate the growth rate of the number of all
possible paths in graphs with respect to their length in an exact
analytical way. Apart from the typical rates of growth, i.e.
exponential or polynomial, we identify conditions for a stretched
exponential type of growth. This is made possible by combining two
or more graphs over the same alphabet, in order to obtain a discrete
dynamical system generated by a triangular map, which can also be
interpreted as a discrete non-autonomous system. Since the vertices
and the edges of a graph usually are used to depict the states and
transitions between states of a discrete dynamical system, the
combination of two (or more) graphs can be interpreted as the
driving, or perturbation, of one system by another.
\\
\\




\end{abstract}


\section{Introduction.}\label{intro}

Strings consisting of sequences of symbols play a central role in
many natural phenomena. For instance, the DNA and RNA are linear
strings carrying the genetic code written in an alphabet consisting
of four letters A, C, G and T (or U) according to whether the
nucleotide subunit contains the base adenine, cytosine, guanine and
thymine (or uracil). Furthermore most of the communication processes
generating and transporting information
 or having a cognitive dimension such as computer programs, the electrical activity of the brain,
texts or music, are implemented in one way or the other in symbolic
terms.

With the advent of nonlinear dynamics and in particular chaos
theory, it has been realized that the evolution of large classes of
dynamical systems can be also described, under certain conditions,
by a sequence of symbols. The existence of such a \emph{symbolic
dynamics} related, in addition, to the principal indicators of the
complexity of the underlying system, provides a unique opportunity
to understand the mathematical and physical mechanisms at the basis
of information generation and transduction. In this context, a
central question is how to enumerate and to characterize the full
set of possible sequences generated by a dynamical system or, in a
more dynamic language, the paths leading from one state to another.

A variety of approaches to this major problem have been reported, but several aspects remain only partially
understood. A measure of the dependence of the number of sequences of a given length (to which we will
refer hereafter as \emph{words}) with respect to the length is the block entropy, defined by

\begin{equation}\label{shannon}
 H(n) = -\sum_{w(n)} P[w(n)] \ln \left (P[w(n)] \right)
\end{equation}

\noindent where the summation is over all allowed words $w(n)$ of
length $n$ and $P[w(n)]$ is the associated probability (a word is
allowed if it represents a possible path leading from one state to
another; in general not all conceivable paths are possible for the
system under investigation). Note that for uniformly distributed
probability, $H(n)$ is exactly the logarithm of the number of words
of length $n$. For sequences generated by Bernoulli or Markov
processes the Shannon-McMillan theorem \cite{Applebaum} asserts that

\begin{equation}\label{sh-mac}
 P[w(n)] \approx  e^{-H(n)} = e^{-n h}
\end{equation}

\noindent where $h$, the \emph{entropy of the source}, is a discrete
analog of the Kolmogorov-Sinai entropy. This entails that (a) $H(n)$
scales linearly with the word length and (b) that long words are
extremely improbable, as they are exponentially penalized. The
penalization is maximal for a Bernoulli processor where $h$ takes
its largest value for a given alphabet \cite{Applebaum},
\cite{Kitchens}. Although still exponential, it is milder for a
Markov processor since $h$ is in this case less than the maximum. As
it turns out holds true for sequences generated by chaotic dynamical
systems as well in the domain of fully developed chaos
\cite{Gaspard}, \cite{NicolisBook}, \cite{Bai-lin90:0}.

In many real world systems we do not
observe exponential rates of growth. This means that there exists a
procedure of \emph{selection} of a subset of sequences possessing
some prescribed properties out of the total number of sequences. In
this context Ebeling and Nicolis [1992] proposed the
following generic scaling law for the block entropy
\begin{equation} \label{ebeling}
H(n)= nh + g n^{\mu}\left( \log n \right)^{\nu}+ e, \quad
h \ge 0,\ 0< \mu <1,\ \nu \le 0
\end{equation}
\noindent and showed that classical literature texts or music
obey to Eq.(\ref{ebeling}) with $h=0, \  \nu=0 $ and $\mu < 1$.
Dynamical systems showing weak chaos in the form of intermittency,
and sporadic systems have also been shown to give rise to a
sublinear scaling of $H(n)$, \cite{Gaspard}, \cite{NicolisBook}.
The question has therefore been raised,
under what conditions could the block entropy of a dynamical system
scale sublinearly with $n$, thereby rendering possible the selection
of long sequences with a non-negligible probability. Obviously one
can think of a number of constraints \emph{external} to the dynamics
(for instance, the \emph{meaning} of a string of letters in a
natural language). A natural and interesting problem concerns the
possibility of constructing a dynamical system including in itself a
selection rule capable of producing a sublinear scale. Our principal
goal in this paper is to propose the construction of higher
dimensional dynamical systems obtained by suitable combinations of
simple systems, each of them generated by a single graph.

The idea
of our construction comes from the theory of two-variable
{\it triangular maps}, i.e., maps where the first component depends only
on the first variable. In our particular case, actually, the system generated by
these triangular maps can also be interpreted
as a {\it non-autonomous} discrete dynamical system.
We are able to show that this mechanism can
generate a subclass of scalings of the form of Eq.(\ref{ebeling}),
leading at the level of Eq.(\ref{sh-mac}) to a stretched exponential
dependence of $P[w(n)]$ on $n$.

To the degree that the states of a system can be encoded via a
reasonably small finite alphabet and as long as the alphabet remains
constant, it is natural to identify the states of the system as
vertices of a graph and the edge between two vertices as a
transition between these two states. The connectivity of a graph
reflects, then, the possible transitions between states, and a path
traversing through the vertices of the graph describes a possible
trajectory of the system. In the present work we shall
adopt this representation.

In this representation the upper bound of the block entropy (see Eq.(\ref{shannon}))
computed over all possible Markov measures \cite{Gaspard} gives the
\emph{topological entropy} of the system, which indicates what paths are
possible to be realized 'in principle' independently of
any probabilistic consideration.

More technically, \cite{IntroSymbDyn}, \cite{Kitchens},
the topological (block) entropy of a sequence X
(or "the entropy of the shift X", if the sequence comes from a shift)
is defined as
\[ h_{top} = \lim_{n \rightarrow \infty } \left(\frac{1}{n} \log_{2} | B_{n}(X) | \right) \]
\noindent where $| B_{n}(X) |$ is the number of n-blocks appearing in X.
Obviously, since all the blocks counted cannot be more than all the combinations of the
letters of the alphabet at hand, $| B_{n}(X) | = |A|^{c \cdot n} \le|A|^{n} $.
\noindent where $|A|$ is the cardinal number of the set of letters $A$
which constitute the alphabet in use.
Usually the parameter $c$ is called entropy's "growth rate".
In this formulation the analogue of Eq.(\ref{ebeling}) would
be the presence of a more intricate relationship, e.g. in the form
$| B_{n}(X) |  = |A|^{c(n)} $.

In Sect.\ref{review}, we review some results for the number of
symbol sequences obtained by following all paths over a graph. Then
in Sect.\ref{combination}, we investigate the possibility of
deviating from this standard symbol sequence generation by
combination of graphs which correspond to a kind of perturbation of
the original dynamical system. We report two examples which give a
stretched exponential growth rate for the resulting symbol
sequences. We conclude in Sect.\ref{conclusion} by an
overview of the results and an outlook of further possible
developments.

\section{Symbolic Dynamics of a single graph}\label{review}
In all generality, the technique of associating a dynamical system
to a symbol-sequence producing
system consists of applying a suitable partition
of the phase space of the system (coarse graining) to a finite number of
cells \cite{NicolisBook}, \cite{IntroSymbDyn}
and then to associate the cells of the partition
to the letters of an alphabet. A trajectory of the system produces
then a symbol sequence as it visits the different cells of the partition.
Relevant partitioning of the phase space  reflects the ability of
recording relevant features of the dynamics.

A key point in the foregoing is to associate a directed graph $\mathbf G$  over the alphabet
of the (coarse grained) states
with the transitions from each state to another.
 The adjancency matrix, $M$, of the graph $\mathbf G$  has
elements $m_{i,j}=1$ if the transition from state $i\rightarrow j$
is possible and $m_{i,j}=0$ otherwise.
Providing, in addition, a Markov measure over the states will give the
transition probability  matrix of the system.

Evidently all information about the  possible
transitions is encoded in the adjancency matrix of the graph and the
growth of the block entropy of the resulting symbol sequences will be
given by the first eigenvalue of $M$. Any path of transitions
associated to a succession of states can be envisioned as a path on the graph $\mathbf G$.
We shall be  interested here in the characterization of possible transitions, and
therefore restrict this work on simple connected graphs associated with
dynamical systems, leaving the graphs equipped with a Markov measure
for later work. In order to proceed some definitions and a review
of known results are necessary.

\vskip 2mm Let us  consider an alphabet $A=\{X_1,X_2,\cdots,X_k\}$
and denote by $\mathcal A_k$ the set of all finite {\bf words}
generated by $A$. By $\mathcal A_k^n$ we denote the set of all words
of length $n$, so we have
\[\mathcal A_k=\cup_{n=1}^{\infty} \mathcal A_k^n.\]

Let $\mathbf G$ be a \emph{ connected directed graph}
(from now on all graphs here considered are connected and directed)
having $A$ as set of
vertices and with at most two edges with different directions
connecting two vertices; thus, there is at most one arrow from a
vertex $X_i$ to another vertex $X_j$. Given a word $s \in \mathcal
A_k$, a $1$--letter extension of $s$ generated by $\mathbf G$ is
defined as follows: if $X_i$ is the last letter of $s$, we add a
letter $X_j$ if and only if the graph $G$ as an arrow from $X_i$ to
$X_j$. Depending on the number of arrows in $\mathbf G$ starting
from $X_i$ we have the same number of $1$--letter extensions of $s$.

As a first step we formulate the generation of words by a
graph in terms of a {\it discrete dynamical system}, i.e., the iterated
application of a suitably defined function over a certain set.
To do this, we define a map $\mathcal G$ from the family $2^{\mathcal A_k}$
of all subset of $\mathcal A_k$ into itself in the following way:
\vskip 2mm {\it for $S \subset \mathcal A_k$, $\mathcal G(S)$ is the
set of all $1$--letter extensions of the words in $S$, generated by
$\mathbf G$.} \vskip 2mm

The set $\mathcal G^{n-1}(A)$ ($\mathcal G^{n-1}$ is the
$n-1$--iterate of $\mathcal G$) is the \emph{ set of the words of
length $n$ generated by $\mathbf G$}. In the following we will
denote this set by $\mathcal W_G^n$ and by $\mathcal W_G$ the set of
all finite words generated by $\mathbf G$ or {\bf admissible words}.
Here, by an admissible word we mean a finite string constructed with
the alphabet $\{X_1,X_2,\cdots,X_k\}$ such that if a pair $X_iX_j$
appears in the string, then $\mathbf G$ has an arrow from $X_i$ to
$X_j$.

With the symbol $\mathcal W_G^n(X_i,X_j)$ we denote the set of
admissible words of length $n$ starting with $X_i$ and ending with
$X_j$. Moreover, we set
\[\mathcal W_G^n(X_i,\cdot)=\cup_{j=1}^k \mathcal W_G^n(X_i,X_j),
\qquad \mathcal W_G^n(\cdot,X_j)=\cup_{i=1}^k \mathcal
W_G^n(X_i,X_j).\]

\vskip 5mm

The following theorem is well known \cite{IntroSymbDyn}:

{\bf Theorem 1} {\it Let $M^{n-1}=[m_{ij}^{(n-1)}]$ be the $(n-1)$--th
power of the adjacency matrix $M$.

The
cardinality  $\omega_G^n(X_i,X_j)$, that is
the number of words in $\mathcal W_G^n(X_i,X_j)$, is then

\[\omega_G^n(X_i,X_j)=m_{ij}^{(n-1)}\]
and
\[\omega_G^n=\sum_{i,j=1}^k m_{ij}^{(n-1)}.\]}

Next we can state a theorem which provides a linear difference equation for
$\omega_G^n$.

{\bf Theorem 2} {\it Let $P_M(\lambda)=\lambda^k+\sum_{r=0}^{k-1}
a_r\lambda^r$ be the characteristic polynomial of $M$
($P_M(\lambda)=det(\lambda I-M)$). Then, for every $n>k$, we have
\[\omega_G^n=-\sum_{r=0}^{k-1} a_r\omega_G^{n-k+r}.\]}

{\bf Proof-} Since any matrix is a root of its characteristic
polynomial, we have
\[M^k=-\sum_{r=0}^{k-1} a_rM^r.\]
Multiplying both sides by $M^{n-k-1}$ we get
\[M^{n-1}=-\sum_{r=0}^{k-1} a_rM^{n-k+r-1}.\]
Thus, we immediately have
\begin{equation}\label{eq1}
\omega_G^n=-\sum_{r=0}^{k-1} a_r\omega_G^{n-k+r}
\end{equation}
{\bf QED}.

\vskip 4mm As an obvious consequence of Theorem 2 we obtain the
following corollary: \vskip 1mm

{\bf Corollary 3} {\it Assume that $\mu_1,\cdots,\mu_s$ are the non
zero distinct roots of the characteristic polynomial of $M$, with
respective multiplicities $r_1,\cdots,r_s$. Then
\[\omega_G^n=c_0+\sum_{h=1}^s \Big(\sum_{q=1}^{r_h}
c_{hq}n^{q-1}\mu_h^n \Big),\]

where the coefficients $c_0$ and $c_{hq}$ are determined by the
values $\omega_G^p$ for\\ $p=1,\cdots,k$.}

\vskip 3mm

We remark that
not only $\omega_G^n$ satisfies the previous difference equation (\ref{eq1}),
but also $\omega_G^n(X_i,X_j)$, for every $i,j=1,\cdots,k$.

\vskip 4mm

We conclude this section by highlighting an important consequence stemming from Corollary 3.
Let us define $\rho=\max \Re(\mu_i) $, then the  rate of growth of
$\omega_G^n$ with respect to the length $n$ of the words can {\it only} have one of the
following three forms:
\begin{enumerate}
\item[(i)]  exponential, $\omega_G^n \simeq \rho^n$, if $\rho >1 $,
\item[(ii)] polynomial, $\omega_G^n \simeq p(n)$,
with $p(n)$ a polynomial with $deg(p(n) \geqslant 1$ , if  $\rho = 1 $ or with $deg(p(n)) =0$
if $\rho = 0 $
\item[(iii)] a combination of polynomials and exponentials, $\omega_G^n \simeq p(n) \rho^n$,
 if $\rho >1 $ and its corresponding root has multiplicity greater than one.
\end{enumerate}

The graphs used in the following section are examples of case (i) and (ii). To get a
polynomial of degree zero, it is enough to consider the graph on a two-letter alphabet
with the following adjacency matrix $M$,
\begin{equation*}
M=
\begin{bmatrix}
0& 1\\
1& 0
\end{bmatrix}
\end{equation*}
\noindent To our knowledge, and according to some preliminary numerical investigation, it seems that
case  (iii) cannot arise.

\noindent The previous remark shows that there is {\it simply no possibility}  for paths
over any graph  to generate symbol-sequences with stretched
exponential growth rates, i.e.,
$\omega_G^n \simeq \rho^{n^{\nu}}$ with $0 < \nu <1$. We shall see,
presently, that this kind of growth may result from a {\it combination}
of graphs.

\section{Combination of graphs.}\label{combination}

Let us now turn to symbol sequences generated by a collection of
graphs. Assume we are given $\ell$ different graphs
$\bf G_1,\bf G_2,\cdots,\bf G_{\ell}$ on the same alphabet
$A=\{X_1,X_2,\cdots,X_k\}$. We also equip this
collection of graphs with a transition rule among them. This rule
might represent, for instance, a 'superselection' rule \cite{JSN2005} not
encoded within each graph but governing the transition among
different domains of distinct dynamics, as in the case of hybrid
systems, see for example \cite{VBThesis} and references therein. We
shall assume that this rule is a deterministic one given by an increasing
sequence of non negative integers, $\{g_i\}$, $i \in \mathbb{N}$,
$g_0=0$.

We construct on $\mathbb{N} \times 2^{\mathcal A_k}$ the following
dynamical system:
\begin{equation}
\begin{split}
\mathcal F(m,S)=(m+1&,\mathcal G_r(S))\quad \mbox{if} \quad
g_{t\ell+r-1}\le m<g_{t\ell+r},\\ &t=0,1,2,\cdots, \quad
r=1,2,\cdots,\ell.
\end{split}
\end{equation}

\noindent This amounts to using repeatedly in sequence the graphs $\bf
G_1,\bf G_2,\cdots,\bf G_{\ell}$, the {\it duration} of the use
of a single graph being controlled by the given deterministic sequence
$\{g_i\}$.

The above construction (analogously to that presented in the previous section
giving a dynamical description of the generation of words by a graph)  gives a
formulation in terms of a discrete dynamical system of the combination of the graphs
 $\bf G_1,\bf G_2,\cdots,\bf G_{\ell}$ driven by the sequence  $\{g_i\}$.
This can be understood as a two-variable map with the first component depending only on the
first variable known as a  {\it triangular map}. The dynamical system is generated by iteration of this map.
Due to the fact that the first variable is discrete,
this construction has also a natural interpretation as a {\it non-autonomous} dynamical system.
In any case, what is important and is in a certain sense the goal of the present investigation,
is to embed the transition rule into the dynamics.

We are mainly concerned with the trajectory starting with $(0,A)$
and, by denoting with $pr_2$ the projection on the second
coordinate, the set

\[\mathcal W_{\mathcal F}^n:=pr_2 (\mathcal F^{n-1}(0,A))\]

\noindent will be called the {\it set of words of length $n$
generated by the dynamical system} (5). In the two following subsections
we shall demonstrate how
non-standard, e.g. stretched exponentials growth rates can result by
such a combination of graphs. First, we will give an example where
both graphs have an absorbing state (for our purposes, an absorbing
state is a vertex with no outcoming edges other
than its own self). Secondly, we give an example of a combination of
a fully connected graph and a graph with an absorbing state.

\subsection{Combination of two graphs: An example with absorbing states}

The two graphs to be combined are defined over the three--letter alphabet
$\{X,Y,Z\}$ as follows.
The first graph $\mathbf G_1$ has the form:

\vskip 10mm
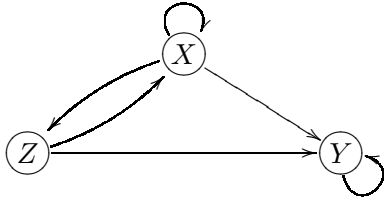
\begin{centering}\hfill
\xymatrix@=40pt@M=5pt@R=20pt{&*+[o][F]{X} \ar@(ul,ur)[] \ar@/_/[dl]
\ar@{->}[dr]
\\*+[o][F]{Z}\ar@{->}[rr] \ar@/_/[ur]
&{} &*+[o][F]{Y} \ar@(d,r)[]\\
}\hfill
\end{centering}
\vskip 5mm

\noindent Its adjacency matrix is

\begin{equation*}
M_1=
\begin{bmatrix}
1& 1& 1\\
0& 1& 0\\
1& 1& 0
\end{bmatrix}
\end{equation*}

\noindent and the corresponding characteristic polynomial is
\[P_{M_1}=\lambda^3-2\lambda^2+1,\]
which has $1$, $(1+\sqrt 5)/2=:\mu$ and $(1-\sqrt 5)/2=:1-\mu$ as
roots.

By applying the results presented in Section 2, we obtain the
following formulas (the subscript $1$ to all $\omega$'s refers to
the graph $\mathbf G_1$):

\begin{equation*}
\begin{aligned}
&\omega_1^n(X,X)=\frac{1}{\sqrt 5}[\mu^n-(1-\mu)^n]\\
&\omega_1^n(X,Y)=-2+\frac{5+3\sqrt 5}{10}\mu^n+\frac{5-3\sqrt 5}{10}(1-\mu)^n\\
&\omega_1^n(X,Z)=\frac{5-\sqrt 5}{10}\mu^n+\frac{5+\sqrt 5}{10}(1-\mu)^n\\
&\omega_1^n(X,\cdot)=-2+\frac{5+2\sqrt 5}{5}\mu^n+\frac{5-2\sqrt 5}{5}(1-\mu)^n\\
\\
&\omega_1^n(Y,X)=\omega_1^n(Y,Z)=0, \quad \omega_1^n(Y,Y)=1\\
&\omega_1^n(Y,\cdot)=1\\
\\
\end{aligned}
\end{equation*}

\begin{equation*}
\begin{aligned}
&\omega_1^n(Z,X)=\frac{5-\sqrt 5}{10}\mu^n+\frac{5+\sqrt 5}{10}(1-\mu)^n\\
&\omega_1^n(Z,Y)=-1+\frac{5+\sqrt 5}{10}\mu^n+\frac{5-\sqrt 5}{10}(1-\mu)^n\\
&\omega_1^n(Z,Z)=\frac{3\sqrt 5-5}{10}\mu^n-\frac{3\sqrt 5+5}{10}(1-\mu)^n\\
&\omega_1^n(Z,\cdot)=-1+\frac{5+3\sqrt 5}{10}\mu^n+\frac{5-3\sqrt 5}{10}(1-\mu)^n\\
\\
&\omega_1^n(\cdot,X)=\frac{5+\sqrt 5}{10}\mu^n+\frac{5-\sqrt 5}{10}(1-\mu)^n\\
&\omega_1^n(\cdot,Y)=-2+\frac{10+4\sqrt 5}{10}\mu^n+\frac{10-4\sqrt 5}{10}(1-\mu)^n\\
&\omega_1^n(\cdot,Z)=\frac{1}{\sqrt 5}[\mu^n-(1-\mu)^n]\\
\\
&\omega_1^n=-2+\frac{15+7\sqrt 5}{10}\mu^n+\frac{15-7\sqrt
5}{10}(1-\mu)^n\\
\end{aligned}
\end{equation*}

\vskip 2mm \noindent The following two inequalities will be useful
for further developments:
\begin{equation}\label{ineq3}
\omega_1^n(X,X)>\omega_1^n(X,Z), \qquad
\omega_1^n(Z,X)>\omega_1^n(Z,Z).
\end{equation}

\vskip 5mm The second graph $\mathbf G_2$ is the following:

\vskip 4mm
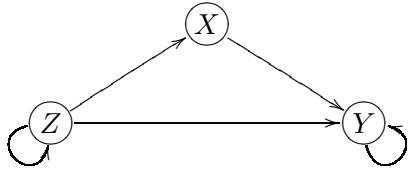
\begin{centering}\hfill
\xymatrix@=40pt@M=5pt@R=20pt{&*+[o][F]{X} \ar@{->}[dr]
\\*+[o][F]{Z}\ar@{->}[rr] \ar@{->}[ur] \ar@(l,d)[]
&{} &*+[o][F]{Y} \ar@(d,r)[]\\
}\hfill
\end{centering}
\vskip 4mm

\noindent Its adjacency matrix is

\begin{equation*}M_2=
\begin{bmatrix}
0& 1& 0\\
0& 1& 0\\
1& 1& 1
\end{bmatrix}
\end{equation*}

\noindent and the corresponding characteristic polynomial is
\[P_{M_2}=\lambda^3-2\lambda^2+\lambda,\]
which has $0$ and $1$ (with multiplicity $2$) as roots.

\noindent We obtain the following formulas:

\begin{equation*}
\begin{aligned}
&\omega_2^n(X,X)=\omega_2^n(X,Z)=0, \quad \omega_2^n(X,Y)=1\\
&\omega_2^n(X,\cdot)=1\\
\\
&\omega_2^n(Y,X)=\omega_2^n(Y,Z)=0, \quad \omega_2^n(Y,Y)=1\\
&\omega_2^n(Y,\cdot)=1\\
\\
&\omega_2^n(Z,X)=\omega_2^n(Z,Z)=1, \quad \omega_2^n(Z,Y)=2n-3\\
&\omega_2^n(Z,\cdot)=2n+1\\
\\
&\omega_2^n(\cdot,X)=\omega_2^n(\cdot,Z)=1, \quad
\omega_2^n(\cdot,Y)=2n-1\\
\\
&\omega_2^n=2n+1\\
\end{aligned}
\end{equation*}

\vskip 4mm
Let now $\{g_i\}$, $i \in \mathbb{N}$, $g_0=0$,  be a
given increasing sequence of non negative integers.
We construct on $\mathbb{N} \times 2^{\mathcal A_3}$ the dynamical
system:
\begin{equation*}
\mathcal F(m,S)=\left\{ \begin{aligned} &(m+1,\mathcal G_1(S)) \quad
\mbox{if} \quad g_{2t}\le m<g_{2t+1}, \quad t=0,1,2,\cdots \\
&(m+1,\mathcal G_2(S)) \quad \mbox{if} \quad g_{2t+1}\le m<g_{2t+2},
\quad t=0,1,2,\cdots
\end{aligned}
\right.
\end{equation*}

\noindent and we denote by $\omega_{\mathcal F}^n$ the cardinality
of the set $\mathcal W_{\mathcal F}^n$. Our aim is to give a lower
and an upper bound for $\omega_{\mathcal F}^n$ and then to prove
that, for suitable choices of the sequence $\{g_i\}$, we may obtain
a growth given by stretched exponentials.

\vskip 3mm

Given two square matrices $M$ and $N$ of the same dimension, we say
that $M \le N$ if this order relation is valid entry by entry. We
write $M<N$ if $M \le N$ and $M \ne N$.

\noindent We denote by $|M|$ the sum of all entries of the matrix
$M$. Clearly, if $M<N$, then $|M|<|N|$. Moreover, if $M<N$, $L<K$
and all entries are non negative, then $ML \le NK$ and $|ML| \le
|NK|$.

\vskip 3mm \noindent If we define $s_t=g_t-g_{t-1}$, we have the
following theorem.

\vskip 2mm

{\bf Theorem 3} {\it For every $t\ge 1$ we have
\begin{equation}\label{eq4}
5^{-t/2}
\prod_{i=1}^t\big[\mu^{s_{2i-1}}-(1-\mu)^{s_{2i-1}}\big]<\omega_{\mathcal
F}^{g_{2t}}<\omega_1^{s_1+s_3+\cdots+s_{2t-1}}\big[1+2\sum_{i=1}^t
s_{2i}\big]
\end{equation}
}

{\bf {Proof}}-- We begin with the second inequality. As a first step
we prove that
\begin{equation}\label{ineq5}
\begin{split}
M_1^{s_1-1}M_2^{s_2}&M_1^{s_3}\cdots
M_2^{s_{2t}}<M_1^{s_1+s_3+\cdots+s_{2t-1}-1}+C(s_1,s_2)+\\
&C(s_1+s_3,s_4)+\cdots+C(s_1+s_3+\cdots+s_{2t-1},s_{2t}),
\end{split}
\end{equation}
where
\begin{equation*}C(a,b)=
\begin{bmatrix}
0& 2b\omega_1^a(X,Z)& 0\\
0& 0& 0\\
0& 2b\omega_1^a(Z,Z)& 0
\end{bmatrix}.
\end{equation*}

\vskip 2mm \noindent The matrices $M_1^r$ and $M_2^r$ have the
following forms: \vskip 2mm
\begin{equation*}M_1^r=
\begin{bmatrix}
\omega_1^{r+1}(X,X)& \omega_1^{r+1}(X,Y)& \omega_1^{r+1}(X,Z)\\
0& 1& 0\\
\omega_1^{r+1}(Z,X)& \omega_1^{r+1}(Z,Y)& \omega_1^{r+1}(Z,Z)
\end{bmatrix}
\end{equation*}
\vskip 2mm

\begin{equation*}M_2^r=
\begin{bmatrix}
0& 1& 0\\
0& 1& 0\\
1& 2r-1& 1
\end{bmatrix}.
\end{equation*}

\vskip 2mm \noindent The proof is by induction on $t$. For $t=1$ we
have:

\begin{equation*}
\begin{split}
M_1^{s_1-1}M_2^{s_2}=&
\begin{bmatrix}
\omega_1^{s_1}(X,X)& \omega_1^{s_1}(X,Y)& \omega_1^{s_1}(X,Z)\\
0& 1& 0\\
\omega_1^{s_1}(Z,Z)& \omega_1^{s_1}(Z,Y)& \omega_1^{s_1}(Z,Z)
\end{bmatrix}+\\
&
\begin{bmatrix}
0& \omega_1^{s_1}(X,X)+(2s_2-1)\omega_1^{s_1}(X,Z)& 0\\
0& 0& 0\\
0& \omega_1^{s_1}(Z,X)+(2s_2-1)\omega_1^{s_1}(Z,Z)& 0
\end{bmatrix}<\\
&
\begin{bmatrix}
\omega_1^{s_1}(X,X)& \omega_1^{s_1}(X,Y)& \omega_1^{s_1}(X,Z)\\
0& 1& 0\\
\omega_1^{s_1}(Z,X)& \omega_1^{s_1}(Z,Y)& \omega_1^{s_1}(Z,Z)
\end{bmatrix}+\\
&
\begin{bmatrix}
0&2s_2\omega_1^{s_1}(X,Z)& 0\\
0& 0& 0\\
0& 2s_2\omega_1^{s_1}(Z,Z)& 0
\end{bmatrix}=M_1^{s_1-1}+C(s_1,s_2).
\end{split}
\end{equation*}

\noindent In order to obtain the previous inequality we used
inequalities (\ref{ineq3}).

\newpage

\noindent Now, assume that the relation is true for $t$ and consider
$t+1$:

\begin{equation}\label{eq6}
\begin{split}
M_1^{s_1-1}&M_2^{s_2}\cdots
M_1^{s{2t}-1}M_2^{s{2t}}M_1^{s_{2t+1}}M_2^{s{2t+2}}<\\
&\big[M_1^{s_1+s_3+\cdots+s_{2t-1}-1}+C(s_1,s_2)+
\cdots+C(s_1+s_3+\cdots+s_{2t-1},s_{2t})\big]\times\\
&\big[M_1^{s_{2t+1}}+C(s_{2t+1}+1,s_{2t+2})\big].
\end{split}
\end{equation}

\noindent Since we have the following equations:

\begin{equation*}
\begin{split}
&C(a,b)C(c,d)=0, \qquad C(a,b)M_1^r=C(a,b),\\
&M_1^{s_1+s_3+\cdots+s_{2t-1}-1}C(s_{2t+1}+1,s_{2t+2})=C(s_1+s_3+\cdots+s_{2t+1},s_{2t+2}),
\end{split}
\end{equation*}

\noindent the last term in Eq.(\ref{eq6}) becomes

\begin{equation*}
M_1^{s_1+s_3+\cdots+s_{2t+1}-1}+C(s_1,s_2)+
\cdots+C(s_1+s_3+\cdots+s_{2t+1},s_{2t+2}).
\end{equation*}

\noindent Thus the inequality (\ref{ineq5}) is proved.

\noindent Since $|C(a,b)|<2b\omega_1^a$, we obtain

\begin{equation*}
\begin{split}
\omega_{\mathcal F}^{g_{2t}}&<\omega_1^{s_1+s_3+\cdots+s_{2t-1}}+2s_1\omega_1^{s_1}+\cdots+2s_{2t}\omega_1^{s_1+s_3+\cdots+s_{2t-1}}\\
\end{split}
\end{equation*}

\noindent Since, obviously,
\begin{equation*}
\omega_1^{s_1+s_3+\cdots+s_{2j-1}}  \leq \omega_1^{s_1+s_3+\cdots+s_{2t-1}},
\end{equation*}
\noindent for any $j \leq t$, we finally get
\begin{equation*}
\begin{split}
\omega_{\mathcal
F}^{g_{2t}}&\omega_1^{s_1+s_3+\cdots+s_{2t-1}} \big[1+2\sum_{i=1}^t s_{2i}\big]
,
\end{split}
\end{equation*}
\noindent i.e., the right-hand side of the above is nothing more than
the right-hand side of Eq.(\ref{eq4}). The other inequality
follows easily, since we have

\begin{equation*}
M_1^{s_1-1}M_2^{s_2}>
\begin{bmatrix}
\omega_1^{s_1}(X,X)& 0& 0\\
0& 0& 0\\
0& 0& 0
\end{bmatrix}
\end{equation*}
\noindent and so
\[|M_1^{s_1-1}M_2^{s_2}M_1^{s_3}\cdots
M_2^{s_{2t}}|>\omega_1^{s_1}(X,X)\omega_1^{s_3+1}(X,X)\cdots\omega_1^{s_{2t-1}+1}(X,X).\]
\noindent Since $\omega_1^{r+1}(X,X) \ge \omega_1^r(X,X)$, we have
\[\omega_1^{s_1}(X,X)\omega_1^{s_3+1}(X,X)\cdots\omega_1^{s_{2t-1}+1}(X,X)\ge
\omega_1^{s_1}(X,X)\omega_1^{s_3}(X,X)\cdots\omega_1^{s_{2t-1}}(X,X)\]
\noindent and
\[\omega_1^{s_1}(X,X)\omega_1^{s_3}(X,X)\cdots\omega_1^{s_{2t-1}}(X,X)=5^{-t/2}
\prod_{i=1}^t\big[\mu^{s_{2i-1}}-(1-\mu)^{s_{2i-1}}\big].\] {\bf
QED}

\vskip 3mm

Now we specify the sequence $\{s_t\}$ (or $\{g_t\}$) in order to
obtain the desired result. Choosing
\[ s_1=4, \quad s_{2t-1}=2t+1, \qquad s_{2t}=(t+1)^4-t^4+t^2-(t+1)^2\]
\noindent we have:
\begin{equation*}
\begin{split}
&s_1+s_3+\cdots+s_{2t-1}=(t+1)^2\\
&s_2+s_4+\cdots+s_{2t}=(t+1)^4-(t+1)^2\\
&n:=g_{2t}=s_1+s_2+\cdots+s_{2t}=(t+1)^4.
\end{split}
\end{equation*}

\noindent Substituting in Eq.(\ref{eq4}) we obtain the inequalities:

\begin{equation*}
5^{-({\sqrt[4]{n}-1})/2}
\prod_{i=1}^{\sqrt[4]{n}-1}\big[\mu^{s_{2i-1}}-(1-\mu)^{s_{2i-1}}\big]<\omega_{\mathcal
F}^n<\omega_1^{\sqrt n}\big[1+2n-2\sqrt n \big].
\end{equation*}

\noindent Denoting by $f_1(n)$ and $f_2(n)$ the left- and
right--hand sides of \ref{ineq7} respectively, we obtain the
following asymptotic evaluations:
\[ f_1(n)\asymp 5^{-({\sqrt[4]{n}-1})/2} \mu^{\sqrt n}, \qquad
f_2(n)\asymp n \mu^{\sqrt n}.\]

\vskip 3mm

At this point it is worthwhile to specify the role of the absorbing
states in the graphs $\mathbf G_1$ and $\mathbf G_2$. In both graphs
the absorbing state if given by $Y$ and the consequence of this is
that $\omega_1^n(Y,\cdot)=\omega_2^n(Y,\cdot)=1$, while
$\omega_1^n(\cdot,Y)=-2+\frac{5+4\sqrt 5}{10}\mu^n+\frac{5-4\sqrt
5}{10}(1-\mu)^n$ and $\omega_2^n(\cdot,Y)=2n-1$. This implies that
when passing from one graph to the other the quite large number of
words ending by $Y$ in not changed along the application of the new
graph. The times of application of each graph are given by $\{ s_{t} \}$ and
in our case this sequence increases sufficiently fast to permit
a strong reduction of the number of generated words.

\subsection{Combination of two graphs: A second example with a fully connected graph}

For this second example we take as graph $\mathbf G_1$ the complete
graph on three letters i.e., the graph:

\vskip 10mm
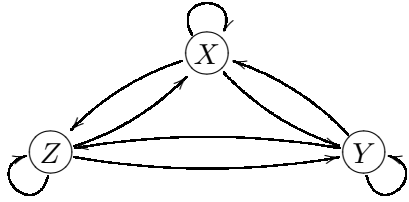
\begin{centering}\hfill
\xymatrix@=40pt@M=5pt@R=20pt{&*+[o][F]{X} \ar@(ul,ur)[] \ar@/_/[dl]
\ar@/_/[dr]
\\*+[o][F]{Z}\ar@(d,l)[] \ar@/_/[ur] \ar@/_/[rr]
&{} &*+[o][F]{Y} \ar@(d,r)[]\ar@/_/[ll] \ar@/_/[ul]\\
}\hfill
\end{centering}
\vskip 5mm

\noindent Clearly this graph has no absorbing states. Its adjacency
matrix is

\begin{equation*}M_1=
\begin{bmatrix}
1& 1& 1\\
1& 1& 1\\
1& 1& 1
\end{bmatrix}
\end{equation*}

\noindent and we immediately have

\begin{equation*}
\begin{aligned}
&\omega_1^n(X,X)=\omega_1^n(X,Y)=\omega_1^n(X,Z)=\\
&\omega_1^n(Y,X)=\omega_1^n(Y,Z)=\omega_1^n(Y,Y)=\\
&\omega_1^n(Z,X)=\omega_1^n(Z,Y)=\omega_1^n(Z,Z)=3^{n-2}\\
&\omega_1^n=3^n\\
\end{aligned}
\end{equation*}

\noindent We remark that $M_1^r=3^{r-1}M_1$.

As second graph we use the graph $\mathbf G_2$ of the previous
example.

We have the following:

{\bf Theorem 4} {\it For every $t\ge 1$ we have
\begin{equation*}
3^{s_1+s_3+\cdots+s_{2t-1}}<\omega_{\mathcal
F}^{g_{2t}}<3^{s_1+s_3+\cdots+s_{2t-1}}\prod_{i=1}^t(2s_{2i}+1).
\end{equation*}
}

{\bf {Proof}}-- We begin by computing $M_1^{s_1-1}M_2^{s_2}$. First
we have the following:

\vskip 2mm
\begin{equation*}M_1M_2^{s_2}=
\begin{bmatrix}
1& 2s_2+1& 1\\
1& 2s_2+1& 1\\
1& 2s_2+1& 1
\end{bmatrix}.
\end{equation*}

\noindent Thus, we obtain the inequalities
\[M_1<M_1M_2^{s_2}<(2s_2+1)M_1.\]

\noindent Hence, we easily get
\[M_1^{s_1-1}M_2^{s_2}>3^{s_1-2}M_1\]
\noindent and, by induction,
\[M_1^{s_1-1}M_2^{s_2}M_1^{s_3}\cdots
M_2^{s_{2t}}>3^{s_1+s_3+\cdots+s_{2t-1}-2}M_1.\]

\noindent So we have
\[|M_1^{s_1-1}M_2^{s_2}M_1^{s_3}\cdots
M_2^{s_{2t}}|=\omega_{\mathcal
F}^{g_{2t}}>3^{s_1+s_3+\cdots+s_{2t-1}}.\]

\vskip 3mm \noindent To prove the other inequality, we see that
\[M_1^{s_1-1}M_2^{s_2}<3^{s_1-2}(2s_2+1)M_1\]

\noindent and, again by induction,
\[M_1^{s_1-1}M_2^{s_2}M_1^{s_3}\cdots
M_2^{s_{2t}}<3^{s_1+s_3+\cdots+s_{2t-1}-2}\Big[\prod_{i=1}^t(2s_{2i}+1)\Big]M_1,\]

\noindent thus

\[M_1^{s_1-1}M_2^{s_2}M_1^{s_3}\cdots
M_2^{s_{2t}}=\omega_{\mathcal
F}^{g_{2t}}<3^{s_1+s_3+\cdots+s_{2t-1}}\Big[\prod_{i=1}^t(2s_{2i}+1)\Big].\]
{\bf QED}

\vskip 4mm

Choosing now the sequence $\{s_t\}$ as in the previous example, we
remark that $s_{2t}=(t+1)^4-t^4+t^2-(t+1)^2 \le 4(t+1)^3$.

\noindent Thus we have:

\begin{equation*}
\begin{split}
\prod_{i=1}^t(2s_{2i}+1)&<\prod_{i=1}^t (3s_{2i})<12^t\prod_{i=1}^t
(i+1)^3=\\
&12^t
\Big(\frac{(t+1)!}{2}\Big)^3=\frac{12^{\sqrt[4]n}}{96}\big[({\sqrt[4]n})!\big]^3.
\end{split}
\end{equation*}

\noindent By using Stirling's formula we obtain for the upper bound
$f_2(n)$ of $\omega_{\mathcal F}^{n}$ the relation
\begin{equation*}
f_2(n)\asymp \sqrt[8]{n^3} 3^{\sqrt n+3/4 \sqrt[4]{n}\log n\log 3}.
\end{equation*}

For the lower bound $f_1(n)$ of $\omega_{\mathcal F}^{n}$ we have
\[f_1(n)=3^{\sqrt n}.\].

\vskip 5mm

The choice of the graphs is obviously crucial for getting
this result. Clearly, by combining graphs each of which produces an
exponential growth rate, we cannot obtain a growth rate given by a
stretched exponential. The two examples presented herein permit to
compute analytically the bounds for the number of words generated by
their combination. In general one can realistically expect only an
estimation obtained by numerical simulation. Furthermore, the choice of the
sequence $\{s_t\}$ has been carefully made, in order to have  the
graph with the polynomial growth rate acting for a longer time. By
combining more than two graphs one has much more freedom in the
choice of the controlling sequence.

\section{Conclusions and Outlook}\label{conclusion}

We have investigated the growth rate of the number of all possible
paths in a graph with respect to their length in an exact analytical
way. By combining two or more graphs over the same alphabet, in order to obtain
a discrete dynamical system generated by a triangular
map, we were able to demonstrate stretched exponential growth of
the resulting symbol-sequences with respect to their length.
The combination of two graphs may be interpreted as the driving of one system by another through
parametric perturbations.
In this view a sub-exponential growth rate indicates 'persistent memory effects'.

In most cases persistent memory makes any finite graph description
only approximate as in the case of zero Kolmogorov-Sinai entropy systems.
The best studied example is that of the Fibonacci system, and the
binary (coarse grained) representation of the circle-map
\cite{Berthe1}, where the 'memory effect' in the language of
dynamical systems excludes any simple graph generating mechanism.
Indeed, symbol sequence of the Fibonacci system can only afford a
description by an {\it infinite}, but countable, family (the leaves
of a rooted tree) of Rauzy graphs, see for example \cite{Casaigne},
 \cite{Dekking:1}, \cite{Dekking:2}. Scaling behavior of
entropy estimates and compressibility issues of such systems with
their well known $ \varpropto log(n) $ behavior, have been
extensively studied \cite{Queffelec}, \cite{Schurmann}, in a
non-probabilistic framework. A combination of graphs might also be
useful in this context. It might account for another description of
such systems with 'long memory'.

The two examples presented here are part of a wide class which
can produce stretched exponential growth. Indeed, it is clear that
by combining more graphs, some of them with polynomial rate of
growth of the number of words, we gain in degrees of freedom in the
choice of the controlling sequence $\{g_t\}$ and thus in lowering
the number of admissible words. It should be noted that in general a
subword of an admissible word generated by a combination of graphs
is not an admissible word.
This contrasts with the case of a single
graph and may be thought as a "mark" of this type of construction.
It is only a matter of direct observation to realize that in natural
languages, too, a subword is not, in general, a word of the language.

Although a great deal of attention has been devoted to the scaling laws
associated with the growth of networks,  attention on
the dynamical basis of the generator of networks from a dynamical systems'
point of view has been reported only quite recently \cite{anselmo2005}.
Now, since the edges of $G$, connecting two states/letters of the
alphabet can be labeled by these states:  i.e. the edge from $i$ to
$j$ has the label $"ij"$, then from the list
of all edges of $G$ we can construct another graph $G^{(2)}$ having
the edges of $G$ as its vertices. This way we can designate the
3-path $"ijk"$ describing the transition $i\rightarrow j \rightarrow
k$, as an edge from $ij$ $jk$. We label this edge from vertex $ij$
of $G^{(2)}$ to the vertex  $jk$ of $G^{(2)}$ as the edge $ij(j)k
\equiv "ijk"$ of $G^{(2)}$. One can continue in this way produce
higher order graphs from the original graph. The study of
the growth of block entropy of the original sequence, i.e. the
proliferation law of words of length $n$ can then be performed by
studying the growth of the associated higher order graphs of $G$.
The same argument goes for the Markov chains associated with $G$ and
its associated measures, if we are interested not in the {\it
possible} transitions but rather in the {\it probable} ones.
Our work opens therefore the perspective that a link between the concepts of block
entropy and network growth. The key observation
here is that by a combination of graphs, the necessary increase of the number of
non admissible words, for example via suitable "super-selection
rules" as in \cite{JSN2005}, is achieved to such a degree that leads
to the establishment of a stretched exponential growth rates. This
sub-exponential growth is not due to a limited sample or system size
but can be thought of as integral part of the topological
description of the system. The presence of specific constraints,
super-selection rules or parameter forcing might play exactly this
role in a more realistic modeling of chaotic or stochastic systems.

\section*{Acknowledgements}
We thank J.S. Nicolis, J.P. Boon and A.
Garcia-Cant\`u for fruitful discussions.
G.L.F. is supported by COFIN--MIUR, 
V.B. thanks ESF--EPEAEK II and particularly the Program PYTHAGORAS II, 
for funding the above work.

\bibliographystyle{alpha}
\bibliography{BaFoNi}

\end{document}